\documentclass[12pt]{amsart}

\usepackage{amsthm}
\usepackage{amsmath}
\usepackage{amssymb}
\usepackage{eucal}
\usepackage{amsfonts}
\usepackage{amsmath}
\usepackage[Lenny]{fncychap}
\usepackage{enumerate}
\usepackage{amssymb}
\usepackage[english]{babel}
\usepackage[breaklinks=true]{hyperref}
\usepackage[all]{xy}
\usepackage{fancybox}
\usepackage{latexsym,mathrsfs,longtable,lscape,dsfont,yfonts}
\usepackage{srcltx}
\usepackage{soul}

\ifx\pdfpagewidth\undefined 
\usepackage[T1]{fontenc}
\else 
\usepackage[OT1]{fontenc} 
\fi 
\usepackage{amsfonts,amssymb,latexsym,longtable,amsmath}

\newtheorem{thm}{Theorem}[section]
\newcommand{\bthm}{\begin{thm}} \newcommand{\ethm}{\end{thm}}
\newtheorem{prop}[thm]{Proposition}
\newcommand{\bprp}{\begin{prop}} \newcommand{\eprp}{\end{prop}}
\newtheorem{fact}[thm]{Fact}
\newcommand{\bfct}{\begin{fact}} \newcommand{\efct}{\end{fact}}
\newtheorem{prob}[thm]{Problem}
\newcommand{\bprb}{\begin{prob}} \newcommand{\eprb}{\end{prob}}
\newtheorem{quest}[thm]{Question}
\newcommand{\bqtn}{\begin{quest}} \newcommand{\eqtn}{\end{quest}}
\newtheorem{lem}[thm]{Lemma}
\newcommand{\blem}{\begin{lem}} \newcommand{\elem}{\end{lem}}
\newtheorem{claim}[thm]{Claim}
\newcommand{\bclm}{\begin{claim}} \newcommand{\eclm}{\end{claim}}
\newtheorem{cor}[thm]{Corollary}
\newcommand{\bcor}{\begin{cor}} \newcommand{\ecor}{\end{cor}}
\newtheorem{conj}[thm]{Conjecture}
\newcommand{\bcnj}{\begin{conj}} \newcommand{\ecnj}{\end{conj}}
\theoremstyle{definition}
\newtheorem{defn}[thm]{Definition}
\newcommand{\bdfn}{\begin{defn}} \newcommand{\edfn}{\end{defn}}
\newtheorem{spec}[thm]{Specializing}
\newcommand{\bspc}{\begin{spec}} \newcommand{\espc}{\end{spec}}
\newtheorem{rem}[thm]{Remark}
\newcommand{\brem}{\begin{rem}} \newcommand{\erem}{\end{rem}}
\newtheorem{cnv}[thm]{Convention}
\newcommand{\bcnv}{\begin{cnv}} \newcommand{\ecnv}{\end{cnv}}
\theoremstyle{remark}
\newtheorem{exam}[thm]{Example}
\newcommand{\bexm}{\begin{exam}} \newcommand{\eexm}{\end{exam}}
\newcommand{\bpf}{\begin{proof}} \newcommand{\epf}{\end{proof}}

\newtheorem*{rep@theorem}{\rep@title}
\newcommand{\newreptheorem}[2]{%
\newenvironment{rep#1}[1]{%
 \def\rep@title{#2 \ref{##1}}%
 \begin{rep@theorem}}%
 {\end{rep@theorem}}
 }
\makeatother

\newtheorem{thmy}{\textbf{Theorem}}

\newreptheorem{theorem}{Theorem}

\newreptheorem{corollary}{Corollary}

\newcommand{\sA} {{\mathcal A}}
\newcommand{\sB} {{\mathcal B}}
\newcommand{\sC} {{\mathcal C}}

\newcommand{\sE} {{\mathcal E}}

\newcommand{\sL} {{\mathcal L}}

\newcommand{\sP} {{\mathcal P}}

\renewcommand{\phi}{\varphi}
\renewcommand{\theta}{\vartheta}

\newcommand{\rf}{\mbox{\large r}}

\newcommand{\topo}{{\mathbf{Top}}}
\newcommand{\el}{{\mathcal{L}}}
\newcommand{\tel}{{\mathbf{Top}\mathcal{L}}}
\newcommand{\stel}{{\mathbf{STop}\mathcal{L}}}

\newcommand{\skp}{\smallskip}
\newcommand{\mkp}{\medskip}
\newcommand{\bkp}{\bigskip}
\def\defi{\buildrel\rm def \over=}

\begin{document}

\title[Reflections in topological algebraic structures]{Reflections in topological algebraic structures}

\author[J. Hern\'andez-Arzusa]{Julio Hern\'andez-Arzusa}
\address{Universidad de Cartagena, Departamento de Matem\'{a}ticas,
Campus de San Pablo, Cartagena, Colombia.}
\email{jhernandeza2@unicartagena.edu.co}


\author[S. Hern\'andez]{Salvador Hern\'andez}
\address{Universitat Jaume I, Departamento de Matem\'{a}ticas,
Campus de Riu Sec, 12071 Castell\'{o}n, Spain.}
\email{hernande@uji.es}

\thanks{ The second listed
author acknowledges partial support by Universitat Jaume I, grant P11B2015-77;
 Generalitat Valenciana,
grant code: PROMETEO/2014/062; and the Spanish Ministerio de Econom\'ia y Competitividad,
grant MTM-2016-77143-P}

\begin{abstract}
Let $\sC$ be an epireflective category of $\topo$ and let $\rf_\sC$\, be the epireflective functor associated with $\sC$.
If $\sA$ denotes a (semi)topological algebraic subcategory of $\topo$, we study when $\rf_\sC\,(\sA)$ is an
epireflective subcategory of $\sA$. We prove that this is always the case for semi-topological structures and
we find some sufficient conditions for topological algebraic structures. We also study when the epireflective functor
preserves products, subspaces and other properties. In particular, we solve an open question about the coincidence
of epireflections proposed by Echi and Lazaar in \cite[Question 1.6]{Echi:MPRIA} and repeated in \cite[Question 1.9]{Echi:TP}.
Finally, we apply our results in different specific topological algebraic structures.
\end{abstract}

\thanks{{\em 2010 Mathematics Subject Classification.} Primary: 54B30, 18B30, 54D10; Secondary: 54H10, 22A30 \\
{\em Key Words and Phrases: Topological algebraic structure, Epireflection, Separation axiom, Mal'tsev spaces} }


\date{\today}
\maketitle \setlength{\baselineskip}{24pt}

\section{Introduction}

In this paper we deal with some applications of epireflective functors in the investigation of topological algebraic structures.
In particular, we are interested in the following question:  Let $\sC$ be an epireflective  subcategory
of $\topo$ (i.e., productive and hereditary, hence containing the 1-point space as the empty product),
and let $\rf_\sC$ be the epireflective functor associated with $\sC$.
If ${\sA}$ denotes a (semi)topological varietal subcategory of  $\topo$ (that is,
a subcategory that is closed under products, subalgebras and homomorphic images),
we study when  $\rf_\sC(\sA)$ is an epireflective subcategory of $\sA$.
From a different viewpoint, this question has attracted the interest of many researchers recently
(cf. \cite{FucZha,Peng,Ravsky:02,Sanchez,T4,T2,T3,XieLiTu})
and, to some extent, our motivation for this research has been to give a unified approach to this topic. First, we recall some
definitions and basic facts.

A full subcategory $\sA$ of a category $\sB$ is \emph{reflective} if the canonical embedding of $\sA$ in $\sB$
has a left adjoint $\rf_\sA : \sB\longrightarrow \sA$ (called \emph{reflection}). Thus for each $\sB$-object $B$ there exists an $\sA$-object $\rf_\sA B$ and
and a $\sB$-morphism $\rf_{(B,\sA)} \colon B\to \rf_\sA B$ such that for each $\sB$-morphism $f\colon B\to A$ to an $\sA$-object $A$,
there exists a unique $\sA$-morphism ${\overline {f}}\colon \rf_\sA B\to A$ such that the following diagram commutes

\[
\xymatrix{ B \ar@{>}[rr]^{\rf_{(B,\sA)}} \ar[dr]_{f} & & \rf_\sA B\ar[dl]^{\overline {f}} \\ & A &}
\]
\mkp

The pair $(\rf_\sA B,\rf_{(B,\sA)})$ is called the \emph{$\sA$-reflection} of $B$ and the morphism $\rf_{(B,\sA)}$ is called \emph{$\sA$-reflection arrow}.
If all $\sA$-reflection arrows are epimorphisms, then the subcategory $\sA$ is said to be \emph{epireflective}.
The functor $\rf\colon \sB \to \sA$, which is called the \emph{reflector}, assigns to each $\sB$-morphism
$f\colon X\longrightarrow Y$, the $\sA$-morphism $\mbox{r}_\sA(f)$ that is determined by the following commutative diagram

$$
\xymatrix{ X\ar[dd]_{\rf_{(X,\sA)}} \ar[rr]^{f}& & Y\ar [dd]^{\rf_{(Y,\sA)}}\\ \\ \rf_\sA X\ar[rr]^{\mbox{r}_\sA\,(f)} & & \rf_\sA Y.}$$
\bkp

\section{Basic facts}

We collect in this section some known facts about epireflective categories that will be used along the paper.
Here, we look at the category $\topo$ of topological spaces and continuous functions. Following Kennison \cite{keni},
by a \emph{topological property} $\sP$, we mean a full subcategory of $\topo$ which is closed under the formation of isomorphic
( = homeomorphic) objects. In general, such subcategories are called \emph{isomorphism-closed},
or \emph{replete}. When speaking about \emph{subcategories}, we will always suppose that they
are full and replete. A topological property $\sP$ is
\emph{hereditary (resp. divisible, productive,
or coproductive)} if the objects of $\sP$ are closed under the formation of subspaces
(resp. quotient spaces, product spaces, or coproduct spaces.)
Here, the terms ``product'', ``coproduct'', ``quotient space'' and ``subspace'' are used in their
topological sense. In particular, Kennison proved that a full subcategory $\sP$ of $\topo$ is epireflective
if and only if $\sP$ is hereditary and productive (cf. \cite{keni}).
Well known examples of reflective subcategories in $\topo$ are:  the classes of all
$T_0$, $T_1$, $T_2$, and $T_3$ spaces, the class of all regular spaces, the completely regular
spaces, the class of all totally disconnected spaces (cf. \cite{keni}). For $T_{3.5}$ spaces
there exist the following reflections: the Stone-\v Cech compactification, the Hewitt realcompactification.
\mkp

Let $\sC$ denote an epireflective subcategory of $\topo$. That is,
for each topological space $X$, there exists an associated topological space $\rf_\sC X\in \mathcal{C}$
 and a surjective continuous function   $\rf_{(X,\sC)}\colon X\longrightarrow \rf_\sC X$ such that
  for every continuous function $f\colon X\longrightarrow Y$, with $Y\in \mathcal{C}$, there exists a continuous function
  $\widetilde{f}\colon \rf_\sC X\longrightarrow Y$ (unique by surjectivity of $\rf_{(X,\sC)}$) such that the following diagram commutes

\[
\xymatrix{ X \ar[rr]^{\rf_{(X,\sC)}} \ar[dr]_{f} & & \rf_\sC X\ar[dl]^{\widetilde{f}} \\ & Y &}
\]
\mkp

For a topological space $X$ and a continuous map $f\colon X\longrightarrow Y$ between topological spaces
we write $\rf_\sC X$ and $\rf_\sC(f)\colon \rf_\sC X\longrightarrow \rf_\sC Y$ for the
image of $X, Y$ and $f$ by the epireflection functor associated to $\sC$. If the subcategory
$\sC$ is clear from the context, we omit it from the index.
Furthermore, the space $ \rf_\sC X$ is uniquely determined up to a homeomorphism, and $\rf_{(X,\sC)}$ is uniquely determined up to composition
from the right with a homeomorphism. It is well known that for every subcategory $\sA$ of $\topo$
there exists a smallest epi-reflective subcategory $\sC(\sA)$ in $\topo$
containing $\sA$, namely spaces homeomorphic to subspaces of products of spaces from $\sA$, (cf. \cite{keni}).
It is said that $\sA$ \emph{generates} $\sC(\sA)$. In case there is a single space $X$ with
$\sC(\{X\}) =\sC$, we say that $X$ \emph{generates} $\sC$ and $\sC$ is called \emph{simply generated} by $X$.
(cf. \cite{herrlich:1969,HS}).
For example the class $\bf{Top_0}$ of $T_0$ spaces is generated by the Sierpi\'{n}ski space. However, the class $\bf{Top_1}$ of $T_1$ spaces
is not simply generated. In fact, the class $\bf{Top_1}$ is generated by the family of cofinite spaces.
Furthermore, given an infinite cardinal $\kappa$, the space $\kappa_{cof}$, which is the set $\kappa$ equipped with
the cofinite topology, simply generates the epireflective subcategory of $\topo$ generated
by all $T_1$-spaces of cardinality at most $\kappa$ (see \cite{Giuli,herrlich:1969}).

For self-completeness, we recall below a realization of the epireflection associated to a epireflective subcategory $\sC(\sA)$
that is generated by a subcategory $\sA$  of $\topo$.

Let $\mathbf{F}(X,\sA)$ denote the class of all continuous functions of $X$ onto spaces in $\mathcal{A}$. We set the following
equivalence relation: for $f\colon X\to Y$ and $g\colon X\to Z$  in $\mathbf{F}(X,\sA)$,
it is said  that $f$ and $g$ are equivalent, $f\sim g$,
if there is a homeomorphism $\psi\colon Y\to Z$ such that $g=\psi\circ f$.
Set $\kappa=|X|$. Since every continuous image of $X$ can be considered as a
subset of $\kappa$, the family of equivalence classes $\mathbf{E}(X,\sA)=\mathbf{F}(X,\sA)/{\sim}$ defines a set.
Let $E=E(X,\sA)$ be the set defined by selecting a fixed element in each equivalence class in $\mathbf{E}(X,\sA)$ and
let $\varphi_{_{\mathcal A(X)}}=\Delta_E X\colon X\to \prod_{f\in E} f(X)$ be the diagonal map of $X$ into the product
$\Pi_E X=\prod_{f\in E} f(X)$.
We have that $\varphi_{_{\mathcal A(X)}}$ is a continuous function from $X$ into
$\Pi_E X$ and, since $\Pi_E X\in\sE(\sA)$, it follows that $\varphi_{_{\mathcal A(X)}}(X)\in\sE(\sA)$.
It is easy to check that $(\mathcal{A}(X),\varphi_{_{\mathcal A(X)}})$ satisfies the universal property of a reflection.
Indeed, let $h\colon X\longrightarrow Y$ be a continuous function from $X$ into $Y\in\sA$.
Then, there exists $f\in\mathbf{F}(X,\sA)$, say $f\colon X\to Z$, such that $f\sim h$.
Let $\psi\colon Z\to Y$ be a homeomorphism with $h = \psi\circ f$ and let $\pi_f$ be the canonical projection of $\Pi_E X$ in $f(X)$.
We have $f=\pi_f\circ \varphi_{_{\mathcal A(X)}}$, which yields $h=\psi\circ(\pi_f\circ \varphi_{_{\mathcal A(X)}})=(\psi\circ\pi_f)\circ \varphi_{_{\mathcal A(X)}}$.
The general case, when $Y\in\sE(\sA)$, follows easily observing that $Y$ is a subspace of a product of spaces in $\sA$.
\mkp

The following facts are easily verified.

\bprp
Let $\mathcal{C}_{1}$ and $\mathcal{C}_{2}$ be epireflective subcategories in $\topo$ such that $\mathcal{C}_{2}\subseteq\mathcal{C}_{1}$. Then
the pair 

$(\rf_{\sC_2}(\rf_{\sC_1}X),\rf_{(\rf_{\sC_1}X,\sC_2)})$ is a realization of the $\sC_2$-reflection
of $X$ in $\topo$.
\eprp
\mkp

\bdfn
A class $\sC$ in $\topo$ is \emph{closed under supertopologies} if whenever $(X,\tau)\in \sC$ and
$\rho$ is a topology on $X$ finer than $\tau$, it follows that $(X,\rho)\in \sC$.
\edfn
\mkp

The following result, whose proof is folklore, clarifies the action the epi-reflection functor for subcategories closed under supertopologies.
(see \cite{HS} for the proof, which is straightforward anyway).

\bthm\label{Th_supertopolgies}
An epireflective subcategory $\sC$ in $\topo$ is closed under supertopologies
if and only if the reflection arrow $\rf_{(X,\sC)}$ is a quotient mapping.
\ethm
\mkp

A topological space is called \emph{functionally Hausdorff}, or \emph{Urysohn},
if distinct points can be separated by a real-valued continuous function, or if any
two distinct points have disjoint closed neighbourhoods, respectively.

\bcor\label{Cor_supertopologies}
The reflection arrow $\rf_{(X,\sC)}$ is a quotient mapping for each of the following subcategories of $\topo$ defined by
the separation axioms: $T_0$, $T_1$, $T_2$, functionally Hausdorff, and Urysohn.
\ecor
\mkp

The next result gives a general realization of the reflection functor for categories whose reflection arrows are quotients.
We omit its easy proof here.

\bprp\label{Pr_Realization}
Let $X$ be a topological space and let $\sC$ denote an epireflective subcategory of $\topo$ whose reflection arrows
are quotient maps. If $R_\sC$ is the intersection of all equivalence relations $R\subseteq X^2$ on $X$ such that $X/R\in\sC$,
then $\rf_\sC X=X/R_\sC$.
\eprp
\mkp

By Theorem \ref{Th_supertopolgies}, the proposition above applies to the epireflections defined by the
separation axioms: $T_0$, $T_1$, $T_2$, functionally Hausdorff, and Urysohn.
In particular, if $\sC_1$ and $\sC_2$ denote the subcategories defined by $T_1$ and $T_2$, we have the following characterization,
whose proof is folklore.

\bprp\label{Pr_T1T2}
Let $X$ be a topological space  and let $R$ be an equivalence relation on $X$. The following assertions are fulfilled:
\begin{enumerate}
\item $X/R$ is $T_1$ if and only if each equivalence class in $R$ is closed in $X$.
\item If the space $X/R$ is Hausdorff then $R$ is closed in $X\times X$.
Conversely, if $R$ is a closed subset of $X\times X$ and in addition the quotient map $\varphi\colon X\to X/R$ is open,
then $X/R$ is Hausdorff.
\end{enumerate}
\eprp

\bcor\label{Co_T1T2}
Let $X$ be a topological space. Then
\begin{enumerate}
\item $\rf_{\sC_1}X=X/R_{\sC_1}$ where $R_{\sC_1}$ is the intersection of all equivalence relations whose equivalence
classes are closed in $X$.
\item Let $R_X$ be the smallest equivalence relation that is closed in $X\times X$. Then $\rf_{\sC_2}X$ is canonically
homeomorphic to $\rf_{\sC_2}(X/R_X)$.
\end{enumerate}
\ecor
\bpf
The verification of (1) is clear.

As for the proof of (2), 
by  Proposition \ref{Pr_Realization}, we know that there is an equivalence relation $R_2$ on $X$ such that
$\rf_{\sC_2}X=X/R_2$. Furthermore, Proposition \ref{Pr_T1T2}(2) implies that $R_2$ must be closed in $X\times X$.
Hence $R_X\subseteq R_2$. Thus we have the following commutative diagram

$$
\xymatrix{ X\ar[dd]_{\rf_{(X,\sC_2)}} \ar[rr]^{\pi_X}& & X/R_X\ar [dd]^{\rf_{(X/R_X,\sC_2)}}\\ \\
\rf_{\sC_2}X\ar@<1ex>[rr]^{\rf_{\sC_2}(\pi_X)}  & & \rf_{\sC_2}(X/R_X) \ar@<1ex>[ll]^{\phi} } $$
\mkp

Here, the map $\phi$ is continuous because $\pi_X$ and $\rf_{X/R_X,\sC_2}$ are both quotient morphisms.
Thus, $\rf_{\sC_2}(\pi_X)$ and $\phi$ are each other inverse, which completes the proof.
\epf

\section{Epireflective categories in topological and semitopological algebraic structures}

So far, only categories of topological spaces have been considered. However, our main interest lies on
topological algebraic categories. Taking the terminology of Hart and Kunen \cite{HK:Fundamenta},
in what follows an \emph{algebraic system} $\el$ is a set (possibly empty or infinite) of symbols
of constants, symbols of functions (every function symbol has a finiteaa
arity $\geq 1$, i.e., $\Phi\colon X^n\to X$) and a set of equations $\Sigma_\el$ that the
elements, functions and constants must satisfy. (E.g., $\forall x\ ex = x$; this contains a
constant sign, a binary operation sign and a variable.) An algebraic system is also
called a \emph{variety}, but further we will use the terminology \emph{algebraic system}.

For a categorical generalization of algebraic, and varietal 
also called monadic  functors, see \cite{AHS,HS2}

A \emph{structure} $\frak U$ for $\el$ is a
set $A$ (the domain) together with elements $c_\frak U$ (of) and functions $\Phi_\frak U\colon A^n\to A$, for $n\geq 1$ a natural
number, the arity of the respective operation, corresponding to the constants and
operations in $\el$ (these we call the \emph{specifications of the constants and operations from $\el$ in $\frak U$}), that satisfy
the equations established in $\Sigma_\el$.
E.g., when we talk about groups, it is understood that
$\el=\{\cdot ,i,1, \Sigma_\el\}$ (symbols of the product, inverse
element, identity and $\Sigma_\el$ denotes the equations that define a group). In general, groups (and other algebraic systems)
are displayed as $\frak U=(A;\cdot ,i,1)$, avoiding the use of the corresponding set of equations $\Sigma_\el$ for short.
Here, only algebraic systems that are specified by a set of equations are considered (cf. \cite{Gratzer}).

For two structures $\frak U$ and $\frak V$ for $\el$ we say that
$f:\frak U\rightarrow \frak V$ \emph{is algebraically an
$\el$-homomorphism} if $f(c_\frak U) = c_\frak B$ for each constant symbol $c$ of $\el$, and
$f(\Phi_\frak U(x_1, . . . ,x_n)) = \Phi_\frak V(f(x_1), . . . , f(x_n))$ for each function symbol $\Phi$ of $\el$, of arity $n\geq 1$.


A \emph{topological structure} (resp. \emph{semitopological structure}) for $\el$ is a pair
$(\frak U,\tau )$ where $\frak U$ is a structure for
$\el$, and $\tau $ is a topology on $A$ making all
functions in $\frak U$ continuous (resp. separately continuous). We write $\frak U$ for
$(\frak U,\tau )$ if the topology is understood.


Let $\frak U$ and $\frak V$ be two (semi)topological structures
of $\el$, and $f:\frak U\rightarrow \frak V$ . The map $f$ is an \emph{$\mathcal L$-homomorphism} from $\frak U$ to $\frak V$ iff $f$ is
continuous and is algebraically an $\el$-homomorphism.


The class consisting of $\el$-topological (resp. $\el$-semitopological) structures
and $\el$-morphisms defines a subcategory of $\topo$ that will be denoted by $\tel$ (resp. $\stel$).
For example, the category of topological groups $\mathbf{TopGrp}$ is
specified by $\sL=(\cdot,i,1)$ with arities $(2,1,0)$.
\mkp

We thank the referee for the information and references that follows next.

\brem\label{referee}
The notions we are dealing with here appeared at the end of the 19-th century
in a text by Whitehead \cite{Whitehead}, but the theory began to develop well only after the theory of
lattices was sufficiently well developed, in the thirties of the 20th century. This
is the subject of \emph{universal algebra}. Thus, what we have named an \emph{algebraic system}, is also
called a \emph{variety}. As a general references to universal algebras, it is pertinent to mention here
the volumes \cite{Cohn,Gratzer,Kurosh}.
\erem

\bdfn
Let $\Phi$ be a $n$-ary function on $X$, i.e., $\Phi\colon X^n\to X$, where $n\geq 1$ is an
integer. A  $\Phi$-congruence in $X$ is an equivalence relation $R$ in $X$ such that
if $x_i,y_i\in X$, $i=1,\dots, n$ and $(x_i,y_i)\in R$\, for $i=1,\dots n$, then
$(\Phi(x_1,\dots, x_n),\Phi(y_1,\dots, y_n))\in R$. (Observe that the analogue of this for
constants $c$ this is automatically satisfied, since $(c, c)\in R$).

Let $\mathcal L$ be an algebraic system and let $\frak U$ be a structure for $\mathcal L$.
If $X$ is the domain of $\frak U$ and $R$ is an equivalence relation on $X$ that
is a $\Phi_\frak U$-congruence for all function symbols $\Phi\in \mathcal L$, then
we say that $R$ is an \emph{$\el$-congruence}.
\edfn

The following proposition are well known, we include them for the reader's sake.

\bprp\label{congruence}\cite[Lem. 2, p. 36]{Gratzer}
Let $\Phi$ be a $n$-ary function on $X$ and let $R$ be a $\Phi$-congruence. If
$\pi\colon X\to X/R$ is the quotient map, then there is an $n$-ary map
$\Phi_R\colon (X/R)^n\to X/R$ defined by
$\Phi_R(\pi(x_1),\dots, \pi(x_n))=\pi(\Phi(x_1,\dots, x_n))$.
\eprp
\mkp

\bcor\label{quotient}\cite[Th. 2, p. 58]{Gratzer}
Let $\mathcal L$ be an algebraic system and let $\frak U$ be a structure for $\mathcal L$.
If $X$ is the domain of $\frak U$ and $R$ is an $\el$-congruence on $X$, then $X/R$ is the domain for
a structure $\frak V=\frak U/R$ for $\mathcal L$, with constant symbols $c/R$ for constant
symbols $c\in X$, and with function symbols $\Phi_\frak V$ sa\-tis\-fying
$\Phi_\frak V(\pi(x_1), . . . , \pi(x_n))=\pi(\Phi_\frak U(x_1, . . . ,x_n))$.
Here $\pi\colon X\to X/R$ is the natural quotient map.
\ecor
\mkp

\bprp\label{isomorphism}\cite[Th. 1, p. 57]{Gratzer}
Let $\Phi$ and $\Psi$ be $n$-ary maps on $X$ and $Y$, respectively. If $f\colon X\to Y$ is a map
such that $f(\Phi(x_1,\dots, x_n))=\Psi(f(x_1),\dots, f(x_n))$ for all $(x_1,\dots, x_n)\in X^n$,
then there is a $\Phi$-congruence $R$ on $X$ and an injective map $\widetilde{f}\colon X/R\to Y$  that makes
the following diagram commutative

\[
\xymatrix{ X \ar@{>}[rr]^{f} \ar[dr]_{\pi} & & Y \\ & X/R\ar[ur]^{\widetilde{f}} &}
\]

\noindent In particular\, $\widetilde{f}(\Phi_R(\pi(x_1),\dots, \pi(x_n)))=\Psi(f(x_1),\dots, f(x_n))$.
\eprp
\mkp

The following result is a generalization of the first isomorphism theorem for arbitrary $\mathcal L$-structures.

\bthm\label{Th_isomorphism}\cite[Th. 1, p. 57]{Gratzer}
Let $f\colon X\to Y$ be an $\mathcal L$-homomorphism from  $\frak U$ to $\frak V$.
Then there is an $\mathcal L$-congruence $R$ on $X$ and an injective map  $\widetilde{f}\colon X/R\to Y$  that makes
the following diagram commutative

\[
\xymatrix{ X \ar@{>}[rr]^{f} \ar[dr]_{\pi} & & Y \\ & X/R\ar[ur]^{\widetilde{f}} &}
\]
\mkp

\noindent where $\widetilde{f}$ is an $\mathcal L$-homomorphism from $\frak U/R$ to $\frak V$.
\ethm
\mkp

The next proposition is essential in many subsequent results.  Here we present a stronger version than our initial result, for finitely many spaces,
which has been kindly offered to us by the referee. Recall that given a family of topological spaces  $\{X_i\}_{i\in I}$,
a map $f\colon X_1\times\dots,\times X_n\longrightarrow Y$ is named \emph{separately continuous} when for every $i\in I$ and
$\vec{x}\in \prod_{j\not= i} X_j$, the map $f_{\vec x}\colon X_i\to Y$, defined by
$f_{\vec x}(x)\defi f(x;\vec x)$, is continuous for every $\vec{x}\in \prod_{j\not= i} X_j$, where the symbol ; is used to mean
that the variable x is placed at the coordinate i.

\bprp\label{44}
Let $\{X_i\}_{i\in I}$ a family of topological spaces and let
$f\colon \prod_{i\in I} X_i\to Y$ be a separately continuous map into a topological space $Y$.
If $\sC$ is an epireflective  subcategory of $\topo$, then there is a (necessarily unique) separately continuous map
$\overline f\colon \prod_{i\in I} \rf_\mathcal{C}X_i\to \rf_\mathcal{C}Y$. 
\eprp
\bpf
Since colimits are preserved by reflections, and actually by all functors
having a left adjoint (cf. Ad\'amek-Herrlich-Strecker, Prop. 18.10 (given there in
dual formulation), or Herrlich-Strecker, Theorem 27.7), it will suffice to prove that
the topology of separate continuity on a product (the so called \emph{cross topology})
is a final structure, for certain maps from the factors of the product,
hence is a colimit of a diagram having objects among the factors of
the product.

Indeed, let $X_i$, for $i\in I$, be non-empty topological spaces (for some factor the empty space the
product of the underlying sets is empty, so there is a unique topology on it.) Here for
simplicity later we omit the index set $I$. We define on the underlying set of $\prod X_i$ the
following topology. Let $\langle x_i\rangle\in\prod X_i$ have the following neighbourhood base.
We take any open sets $G_i$ of $X_i$ containing $x_i$, and then the neighbourhood base of $\langle x_i\rangle$
will consist of the sets $ \prod_{j\not= i} \{x_j\}\times G_i$ (a \emph{cross}).
This is the intersection of the topologies $X_i\times (\prod_{j\not= i} X_j)_{\text{discr}}$,
where $i\in $I is arbitrary, and the second factor has the discrete topology.
Observe that this second factor can be given also in another way. We take the topological sum of
$|\prod_{j\not= i} X_j|$ copies of $X_i$. Observe that the
underlying set functor $U : \colon\topo\to \textbf{Set}$ has both a left adjoint (discrete topologies
on sets) and a right adjoint (indiscrete topologies on sets). Therefore $U$ preserves
all limits and colimits from $\topo$ to $\textbf{Set}$, in particular products.

In order to take the intersection of topologies, which is once more a colimit, we
have to identify each set $\rf_\sC X_i\times \prod_{j\not= i} \rf_\sC X_j$ canonically to the set $\prod_{j\in I} \rf_\sC Xj$.
However, $U$ preserves products, hence this makes the desired identification (by \emph{iteration of
products}).

In particular, the sum diagram, as a colimit, is preserved, therefore we have
the sum of $|\prod_{j\not= i} X_j|$ many copies of $\rf_\sC X_j$ , with the canonical injections from the
$\rf_\sC X_j$s, where $j\not= i$. This is once more the product of $\rf_\sC X_i$ and a discrete space of
cardinality $|\prod_{j\not= i} X_j|$. We ought to reduce this cardinality to $|\prod_{j\not= i} \rf_\sC X_j|$.

Observe that coequalizers, being colimits, are preserved by reflections. However,
the image will be only a coequalizer, but not a cointersection of extremal epimorphisms.
It will be a cointersection only if the domains of the morphisms (to be
coequalized) are the same. This brings in once more that the cardinality  $|\prod_{j\not= i} X_j|$
ought to be reduced to $|\prod_{j\not= i} \rf_\sC X_j|$. Now recall that cointersections can be made
iteratively. Namely, we divide the maps to be coequalized to (equivalence) classes,
take the coequalizer of each class, which is an extremal epimorphism. Now we have
a lot of extremal epimorphisms (as many as there are classes) and we have to take
the cointersection of them, which is once more an extremal epimorphism. This
gives the coequalizer of all considered morphisms.

Now the classes will be formed as follows. One class is formed by all $\langle x_j\, |\, j\not= i\rangle$,
for which the images $\langle \rf_{\sC,X_j}(x_j)\, |\, j\not= i\rangle$ are a fixed point of $\prod_{j\not= i} \rf_\sC X_j$.
\epf

The next corollary follows from Proposition \ref{44}. However, we include the
proof here for the reader's sake.

\bcor\label{4}
Let $\{X_1,\dots , X_n\}$ be finitely many topological spaces and let
$f\colon X_1\times \dots \times X_n\to Y$ be a separately continuous map into a topological space $Y$.
If $\sC$ is an epireflective  subcategory of $\topo$, then there is a (necessarily unique) separately continuous map
$\overline f\colon \rf_\mathcal{C}X_1\times\dots \rf_\mathcal{C}X_n\to \rf_\mathcal{C}Y$ such that
$\overline f((\rf_{(X_1,\sC)}(x_1),\dots ,\rf_{(X_n,\sC)}(x_n)))=\rf_{(Y,\sC)}(f((x_1,\dots , x_n)))$.
\ecor
\bpf
In order to simplify the notation, we treat the case $n = 2$ only, as this is representative for the general case.
The proof for $n > 2$ is obtained proceeding by induction.

For each fixed point $c\in X_2$, the map $f_c\colon X_1\to Y$, defined by $f_c(x_1):=f(x_1,c)$,
is continuous. Accordingly, there exists a continuous function
$\rf_{\sC}\, (f_c)\colon \rf_\sC X_1\to \rf_\sC Y$ that makes commutative the following digram
$$
\xymatrix{ X\ar[dd]_{\rf_{(X,\sC)}} \ar[rr]^{f_c}& & Y\ar [dd]^{\rf_{(Y,\sC)}}\\ \\ \rf_\sC X\ar[rr]^{\rf_{\sC}\, (f_c)} & & \rf_\sC Y}$$
\bkp

\noindent Thus, the map $\rf_{\sC}\, (f_c)$ is defined, for each $\rf_\sC X_1\ni x'_1=\rf_{(X_1,\sC)}(x_1)$, by
$$\rf_{\sC}\, (f_c)(x'_1)= \rf_{\sC}\, (f_c)(\rf_{(X_1,\sC)}(x_1))= \rf_{(Y,\sC)}(f_c(x_1))= \rf_{(Y,\sC)}(f(x_1,c)).$$

Furthermore, we have that if $\rf_{(X,\sC)}(x_1)=\rf_{(X,\sC)}(u_1)$ then
$\rf_{\sC}\, (f)(x_1,c)=\rf_{\sC}\, (f)(u_1,c)$. This crucial fact will be widely applied in what follows.

The equality above implies that the map $\tilde f\colon \rf_\sC X_1\times X_2\to Y$, defined by
$$\tilde f((\rf_{(X_1,\sC)}(x_1),x_2)):=\rf_{(Y,\sC)}(f(x_1,x_2)),$$
is well defined and separately continuous.
Repeating the argument for  $\tilde f^d\colon X_2\to Y$, defined by $\tilde f^d(x_2):=\tilde f(d,x_2)$
with $d\in \rf_\sC X_1$, for each $d\in \rf_\sC X_1$ there exists a continuous map
$\rf_{\sC}\, (\tilde f^d)\colon \rf_\sC X_2\to \rf_\sC Y$, satisfying
$$\rf_{\sC}\, (\tilde f^d)(\rf_{(X_2,\sC)}(x_2))=\rf_{(Y,\sC)}(\tilde f^d(x_2)).$$
It is now clear that the map $\overline f\colon  \rf_\sC X_1\times \rf_\sC X_2\longrightarrow  \rf_\sC Y$
defined by $$\overline f(\rf_{(X_1,\sC)}(x_1),\rf_{(X_2,\sC)} (x_2))=\rf_{(Y,\sC)}(f((x_1,x_2)))$$ is
well defined and separately continuous.
\epf

As a consequence of the previous result, it follows that epireflections respect semitopological structures
in the best possible fashion.

\bprp\label{pr_semi}
Let $\stel$ be a category of semitopological structures. If $\sC$ is an epireflective  subcategory of $\topo$,
then $\rf_\sC(\stel)\subseteq \stel$. That is, for each $\el$-structure $(\frak U,X)\in \stel$, where $X$ is the domain of $\frak U$,
we have that $\rf_\sC X$ is the domain of an $\el$-structure $\frak V$ such that $(\frak V,\rf_\sC X)\in \stel$.
Furthermore, the reflection arrow $\rf_{(X,\sC)}\colon X\to \rf_\sC X$ is an $\el$-homomorphism in $\stel$ for each $(\frak U,X)\in \stel$.
\eprp
\bpf
Let $(\frak U,X)$ be a semitopological structure for $\stel$.
We equip $\rf_\sC X$ with the algebraic structure $\frak V$ built by taking the constants $c_{\frak V}:=\rf_{(X,\sC)}(c_\frak U)$
for all constants $c\in\el$ and, if $\Phi\in \el$ is a separately continuous $n$-ary function symbol, then we apply Corollary \ref{4}
in order to define, for $x'_i\in \rf_\sC X$, with $x'_i=\rf_{(X,\sC)}(x_i)$, for $x_i\in X$, as
$$\Phi_{\frak V}((\rf_{(X,\sC)}(x_1),\dots, \rf_{(X,\sC)}(x_n)))=\rf_{(X,\sC)}(\Phi_\frak U((x_1,\dots, x_n))).$$
This definition implies that
$\Phi_{\frak V}\colon (\rf_\mathcal{C}X)^n \to \rf_\mathcal{C}X$ is a well defined separately continuous
$n$-ary function symbol for all $\Phi\in\el$.
Thus $\rf_\sC X$ is equipped with a semitopological $\sL$-structure canonically inherited from the $\sL$-structure in $X$.
Furthermore, it is also clear that $\rf_{(X,\sC)}$ is a continuous $\el$-homomorphism.
\epf
\mkp

We are now in position of establishing the main result in this section.

\bthm\label{Th_semicontinuous}
Let $\stel$ be a category of semitopological structures. If $\sC$ is an epireflective  subcategory of $\topo$,
then $\rf_\sC(\stel)$ is an epireflective subcategory of $\stel$.
\ethm
\bpf
First, observe that, since $\rf_{(X,\sC)}\colon X\to \rf_\sC X$ is an $\el$-homomorphism for all $\el$-structure $(\frak U,X)\in \stel$,
it follows that the equivalence relation $\{(x_1, x_2)\in X^2 : \rf_{(X,\sC)}(x_1) = \rf_{(X,\sC)}(x_2)$
is an $\el$-congruence in $X$ for all $(\frak U,X)\in \stel$.
By Proposition \ref{pr_semi}, we know that $\rf_\sC(\stel)$ is equipped with a semicontinuous $\sL$-structure, where
the reflection arrows $\rf_{(X,\sC)}$ are epimorphisms in $\topo$.  Thus it will suffice to show that $\rf_\sC$ preserves
$\el$-morphisms.

Let $(\frak U,X)$ and $(\frak U',Y)$ be two semitopological structures in $\stel$ and let
$(\frak V,\rf_\sC X)$ and $(\frak V',\rf_\sC Y)$ the structures in $\stel$, with domains $\rf_\sC X$ and $\rf_\sC Y$
canonically associated to the former by Proposition \ref{pr_semi}.

Let $f\colon X\to Y$ be a continuous $\sL$-morphism. If $\Phi\in \el$ is a separately continuous $n$-ary function symbol,
by the commutativity of the diagram
$$
\xymatrix{ X\ar[dd]_{\rf_{(X,\sC)}} \ar[rr]^{f}& & Y\ar [dd]^{\rf_{(Y,\sC)}}\\ \\ \rf_\sC X\ar[rr]^{\rf_{\sC}\, (f)} & & \rf_\sC Y}$$
\bkp

\noindent it follows
\mkp

\begin{equation}
\begin{split}
\rf_{\sC}\, (f)(\Phi_{\frak V}(\rf_{(X,\sC)}(x_1)\cdots, \rf_{(X,\sC)}(x_n)))&=
\rf_{\sC}\, (f)(\rf_{(X,\sC)}(\Phi_{\frak U}((x_1,\cdots, x_n)))) \\
 & = \rf_{(Y,\sC)} (f(\Phi_{\frak U}((x_1,\cdots, x_n))))\\
 &=\rf_{(Y,\sC)} (\Phi_{\frak U'}((f(x_1)\cdots, f(x_n))))\\
 &=\Phi_{\frak V'}((\rf_{(Y,\sC)}(f(x_1)),\cdots, \rf_{(Y,\sC)}(f(x_n)))\\&=\Phi_{\frak V'}((\rf_{\sC}\, (f)(\rf_{(X,\sC)}(x_1))\cdots,
 \rf_{\sC}\, (f)(\rf_{(X,\sC)}(x_n)))).
\end{split}
\end{equation}
\epf
\mkp

Theorem \ref{Th_semicontinuous} allows us to obtain a neat realization of epireflections whose reflection arrows are quotient maps.
Remark that we have shown that the epireflection functor in $\topo$
coincides with the epireflection functor in $\stel$.

\bcor\label{Co_Realization}
Let $\stel$ be a category of semitopological structures. If $\sC$ is an epireflective  subcategory of $\topo$
whose reflection arrows $\rf_{(X,\sC)}$ are quotient maps then $\rf_\sC X=X/R_\sC$ for all $(\frak U,X)\in\stel$,
where $R_\sC$ coincides with the intersection of all $\sL$-congruences $R$ such that $X/R\in \sC$.
\ecor
\bpf
It suffices to apply Proposition \ref{Pr_Realization} and Theorem \ref{Th_semicontinuous}.
 \epf
 \mkp

On the other hand, using Corollary \ref{Co_T1T2}, we obtain

\bcor\label{Co_T1T22}
Let $\stel$ be a category of semitopological structures and let $X$ be a space in $\stel$. Then
$\rf_{\sC_1} X=X/R_{\sC_1}$ where $R_{\sC_1}$ is the intersection of all  equi\-va\-lence relations  such that equivalence
classes are closed in $X$, furthermore $R_{\sC_1}$ is an $\el$-congruence in $X$.
\ecor
\mkp

The following result was established by Tkachenko \cite[Theorem 3.4]{T2} for semitopological groups.
Here, we obtain a variant of his proof as an application of our results.

\bcor\label{Co_T1S}
Let $\el$ denote the algebraic system defined by groups and let $\stel$ be the corresponding category
semitopological groups (group operation is separately continuous).
Then $\rf_{\sC_1}G=G/H$ for all $G\in \stel$, where $H$ is the intersection of
all closed subgroups of $G$.
\ecor
\bpf
First, remark that if $R$ is an $\el$-congruence defined on $G$, then the $R$-equivalence classes are cosets of the normal subgroup
$H_R\defi \{x\in G : (x,e_G)\in R\}$. Thus, using Corollary \ref{Co_T1T22}, we obtain that
$\rf_{\sC_1}G=G/H$, where $H$ is a normal closed subgroup of $G$.
Let us see that $H$ the smallest closed subgroup of $G$. Indeed if $K$ is a
closed subgroup of $G$, then $G/K$ is a $T_1$ space, therefore if $\pi\colon G\to G/K$ is
the respective quotient mapping, we have a group homomorphism
$\rf_{\sC_1}(\pi)\colon G/H\to G/K$ that makes the following diagram
\[
\xymatrix{ G \ar@{>}[rr]^{\pi} \ar[dr]_{\rf_{(G,\sC_1)}} & & G/K \\ & G/H \ar[ur]_{\rf_{\sC_1}(\pi)} &}
\]
\mkp

\noindent commutative. This obviously means that $H\subseteq K$.
\epf

\section{Products}

In this section we deal with the epireflections that preserve products, which is a crucial fact in order to
study the preservation of topological structures.

Let $\sC$ be an epireflective category in $\topo$ and let $\{X_i\}$ a set (resp. finite set) of topological spaces.
Then we have the following commutative diagram
$$
\xymatrix{\prod X_i\ar[dd]_{\rf_{(\prod X_i,\sC)}} \ar[rr]^{\prod \rf_{(X_i,\sC)}}& & \prod \rf_\sC X_i\ar [dd]^{\mbox{id}}\\
\\ \rf_\sC \prod X_i\ar[rr]^{\rf_{\sC}\, (\prod \rf_{(X_i,\sC)})} & & \prod \rf_\sC X_i}$$
\bkp

\noindent Set $\mu_\sC\defi \rf_{\sC}\, (\prod \rf_{(X_i,\sC)})$. We notice that $\mu_\sC$ is defined
uniquely by the condition that the following diagram commutes for every $j$

\[
\xymatrix{ \rf_\sC\prod X_i \ar@{>}[rr]^{\mu_\sC} \ar[dr]_{\rf_\sC(\pi_{_{X_j}})} & & \prod \rf_\sC X_i\ar[dl]^{\pi_{_{\rf_\sC X_j}}} \\ & \rf_\sC X_j &}
\]
\mkp

\noindent here $\pi_{_{\rf_\sC X_j}}$ denotes the respective canonical projection. If $\mu_\sC$ is a homeomorphism onto $\prod \rf_\sC X_i$ for every family of
topological spaces (resp. finite family of topological spaces)
then we say that \emph{$\sC$ preserves products} (resp. \emph{$\sC$ preserves finite products}).

For example, from \cite[Prop. 6.2.1]{Morita}, there follows the following sufficient condition for the preservation of finite products.

\bprp\label{products}
Let $\sC$ be an epireflective category in $\topo$ and let $\{X_i : i\in I\}$ be a set of topological spaces such that
$\rf_{(X,\sC)}$ is open for all $X\in \topo$. Then $\sC$ preserves the product $\prod\limits_{i\in I} X_i$.
\eprp
\bpf
Clearly the map $\mu_\sC$ defined above is onto. To prove that $\mu_\sC$ is injective, it will suffice to repeat \emph{mutatis mutandis} the argument used by
T. Ishii in \cite[Prop. 6.2.1]{Morita}. Therefore, we have shown that $\mu_\sC$ is a bijection. Moreover, since every reflection arrow
$\rf_{(X,\sC)}$ is open, it follows that the reflection arrow  $\prod\limits_{i\in I} \rf_{(X_i,\sC)}$
is also open. Therefore, we have the commutative diagram

\[
\xymatrix{ & \prod X_i\ar@{>}[dl]_{\rf_{(\prod X_i,\sC)}} \ar[dr]^{\prod \rf_{(X_i,\sC)}}&\\
\rf_\sC\prod X_i\ar@{>}[rr]^{\mu_\sC} & & \prod \rf_\sC X_i }
\]
\mkp

Since both $\rf_{(\prod X_i,\sC)}$ and $\prod \rf_{(X_i,\sC)}$ are open and $\mu_\sC$ is continuous, it follows that $\mu_\sC$
is a homeomorphism.
\epf

It is a well-known fact that the $T_0$-reflection in $\mathbf{Top}$ preserves all products but
the $T_1$, $T_2$ and $T_3$ reflections do not preserve all finite products (cf. \cite[\S 1]{HD:87}).

The following result improves Proposition \ref{4} for epireflections that preserve finite products.
Again, we thank the referee for letting us notice that our results can be extended to infinite products.

\bprp\label{5} Let $\{X_i : i\in I\}$ be a set of topological spaces and
let $f\colon \prod\limits_{i\in I}X_i\longrightarrow Y$ be a continuous map
into a topological space $Y$. If $\sC$ is an epireflective  subcategory of $\topo$ such that the epireflection preserves
the product $ \prod\limits_{i\in I}X_i$, then
there is a (necessarily unique) continuous map
 $\overline f\colon  \prod\limits_{i\in I}\rf_\mathcal{C}X_i\longrightarrow \rf_\mathcal{C}Y$
 such that $$\overline f(\langle\rf_{(X_i,\sC)}(x_i)\rangle )=\rf_{(Y,\sC)}(f(\langle x_i\rangle )).$$
\eprp

\bpf

By hypothesis, the canonical homeomorphism
$$\mu_\sC \colon \rf_\mathcal{C}\prod\limits_{i\in I}X_i \to \prod\limits_{i\in I}\rf_\mathcal{C}X_i$$

\noindent satisfies that the inverse mapping defined by
$\mu_\sC^{-1}(\langle\rf_{(X_i,\sC)}(x_i)\rangle )= \rf_{(\prod\limits_{i\in I}X_i,\sC)}(\langle x_i\rangle)$
is well defined (i.e., it depends only on $\rf_{(X_i,\sC)}(x_i)$ and not on the choice of $x_i$),
and is continuous.
Therefore, we have the following commutative  diagram:

$$
\xymatrix{ \rf_\mathcal{C}\prod\limits_{i\in I}X_i \ar [dr]_{\rf_\sC(f)}  \ar [rr]^{\mu_{\sC}} & &
\prod\limits_{i\in I}\rf_\mathcal{C}X_i \ar [dl]^{\overline f} \ar [ll]^{\mu^{-1}_{\sC}} \\ & \rf_\sC Y  \\
\prod\limits_{i\in I}X_i \ar[uu]^{\rf_{(\prod\limits_{i\in I}X_i,\sC)}} \ar[rr]^{f} & & Y\ar[ul]_{\rf_{(Y,\sC)}}}$$
\mkp

Hence $\overline f = \rf_\sC(f)\circ\mu^{-1}_{\sC}$ is a continuous map. The equality stated in Prop. \ref{5}
follows from the commutativity of the lower square in the diagram.
\epf

The following result improves Theorem \ref{Th_semicontinuous} when the reflection functor preserves products.

\bthm\label{Th_continous}
Let $\tel$ be a category of topological structures. If $\sC$ is an epireflective  subcategory of $\topo$ such
that the epireflection preserves products, then $\rf_\sC(\tel)$ is an epireflective subcategory of $\tel$.
\ethm
\bpf
It suffices to observe that $\Phi_{\rf_\sC X}$ is  continuous for every $n$-ary function symbol $\Phi\in \el$
as a consequence of Proposition \ref{5}.
\epf
\mkp



\bcor
Let $\tel$ be a category of topological structures and let $\sC_0$ be  the class of $T_0$,
then $\rf_{\sC_0}(\tel)$ is an epireflective subcategory of $\tel$.
\ecor
\bpf
The classes of $T_0$ spaces is both epireflective and preserve products in $\topo$.
\epf
\mkp

\brem
Let $\mathbf{TopGrp}$ denote the category of topological groups and let $\mathcal{C}_0$ denote the epireflective
subcategory of $T_0$ spaces in $\topo$. Then $\rf_{\mathcal{C}_0}(\mathbf{TopGrp})$ is an epireflective subcategory
of $\mathbf{TopGrp}$. Furthermore, since every $T_0$ topological group is $T_{3.5}$, it follows that
every member in $\rf_{\mathcal{C}_0}(\mathbf{TopGrp})$ is $T_{3.5}$. In other words, for topological groups,
the $T_0$ reflection and the $T_{3.5}$ reflection coincide.
\erem

There are more general structures than groups where the epireflective subcategories that are closed under
finer topologies preserve products. An example of this are Mal'tsev spaces.
\bdfn
A \emph{Mal'tsev operation} on a topological space $X$ is a map $\Phi \colon X^3\rightarrow X$ satisfying the
identity $\Phi(x,x,y) = \Phi(y,x,x) = y$ for all $x,y\in X$. A space is a topological \emph{Mal'tsev} space if it admits
a continuous Mal'tsev operation.
\edfn
\mkp

For example, if $G$ is a topological group, then the map $(x,y,z)\mapsto xy^{-1}z$ is a
Mal'tsev operation on $G$. Hence every topological group is a Mal'tsev space. In like manner,
a \emph{semitopological Mal'tsev} space is a space that admits a separately continuous Mal'tsev operation.

The classes $\bf TopMlt$ (resp. $\textbf{STopMlt}$) of topological (resp. semitopological) Mal'tsev spaces  are algebraic systems
in $\topo$ with the continuous maps that respect these algebraic structures as arrows.

According to Gartside, Reznichenko and Sipacheva \cite{GRS:T&A} topological Mal'tsev spaces were introduced by Uspenskij in \cite{Usp:1989}
and have subsequently been studied by
several authors. In this section we deal with these spaces and our main motivation is to transfer
much of the behavior of topological groups to Mal'tsev
spaces. We will see that most epireflective functors that preserve the topological group structure
also respect the topological Mal'tsev operation.

The following result is attributed to Mal'tsev \cite{Mal54} by Reznichenko and Uspenskij \cite[4.11]{RezUsp98}.
Our formulation is somewhat more general.

\blem\label{open}
Let $(X,\Phi)$ be a semitopological Mal'tsev space and let $R$ be a $\Phi$-congruence in $X$.
Then the quotient map $\pi\colon X\to X/R$ is open.
\elem

\bprp\label{Mopen}
Let $(X,\Phi)$ be a a semitopological Mal'tsev space and let $\sC$ be an epi-reflective class in $\topo$
closed under supertopologies. Then the reflection arrow $\rf_{(X,\sC)}$ is an open map.
\eprp
\bpf
By Theorem \ref{Th_semicontinuous}, we know that $(\rf_\sC X,\rf_\sC(\Phi))$ is  a semitopological Mal'tsev space.
Furthermore, the reflection arrow $\rf_{(X,\sC)}$ is a $\Phi$-homomorphism. By Theorem \ref{Th_isomorphism},
there is a $\Phi$-congruence $R$ such that the following diagram commutes

\[
\xymatrix{ X \ar@{>}[rr]^{\rf_{(X,\sC)}} \ar[dr]_{\pi} & & \rf_\sC X \\ & X/R\ar[ur]_{\widetilde{\rf}_{(X,\sC)}} &}
\]
\mkp

By Theorem \ref{Th_supertopolgies}, is a quotient $\Phi$-homomorphism. Since $\pi$ is open
by Lemma \ref{open} and $\widetilde{\rf}_{(X,\sC)}$ is one-to-one, it follows that
$\rf_{(X,\sC)}$ is open.
\epf
\mkp

From Proposition \ref{products} and Proposition \ref{Mopen}, we obtain.

\bthm\label{Mproducts}
Let $\sC$ be an epireflective category in $\topo$ that is closed under
finer topologies.
If $\{X_i\}$ is a family of semitopological Mal'tsev spaces, then $\rf_\sC$ preserves the
product of $\{X_i\}$.
\ethm

\bcor\label{Mproducts2}
If $\sC$ denotes the epireflective category of $\topo$ defined by any of the following
separation axioms: $T_0$, $T_1$, $T_2$ and functionally Hausdorff, then $\rf_\sC$
preserves arbitrary products of semitopological Mal'tsev spaces.
\ecor
\mkp

The next result is an application of the techniques developed in this paper.
It shows that the modification of a Mal'tsev space by most separation axioms
is a Mal'tsev space.

\bthm\label{Mal'tsev}
If $\sC$ denotes the epireflective category of $\topo$ defined by any of the following
separation axioms: $T_0$, $T_1$, $T_2$ and functionally Hausdorff, and $\mathbf{TopMlt}$ is
the subcategory of topological Mal'tsev spaces. Then $\rf_\sC(\mathbf{TopMlt})$ is an epireflective
subcategory of $\mathbf{TopMlt}$.
\ethm
\bpf
Use Corollary \ref{Mproducts2} and Theorem \ref{Th_continous}.
\epf

\bigskip

\section{Subspaces}

We say that an epireflection \emph{preserves} subspaces if the following holds.
For $i\colon A\hookrightarrow X$ an inclusion of a subspace we have
that $\rf_\sC(i)\colon \rf_\sC A \to \rf_\sC X$ is an inclusion of a subspace as well. More exactly: there
is an inclusion of a subspace $i\colon B\hookrightarrow \rf_\sC X$ and and homeomorphism $h$ such that
$j=\rf_\sC(i)\circ h$. In this section, we study when epireflection functors preserve subspaces in topological algebraic structures.

\bdfn
Let $X$ be a topological space. Two subsets $A$ and $B$ of $X$ are said to be {\it completely
separated} in $X$ if there exists a continuous real valued map $f$ such that $f(a) = 0$ for all $a\in A$
and $f(b) = 1$ for all $b\in B$.

The space $X$ is said to be \emph{functionally Hausdorff} if any two different
points can be completely separated by a real valued continuous function.

The space $X$ is said to be {\it completely regular} if every closed set $F$ of $X$ is completely
separated from any point $x\notin F$. A completely regular $T_1$-space is called
a $T_{3,5}$-space.

A topological space $X$ is said to be \emph{Urysohn} if any
two distinct points have disjoint closed neighhbourhoods.
\edfn
\mkp

It is easy to show that epireflection functors do not preserve subspaces in general. In order to see this, consider
the epireflective subcategories defined by the separation axioms $T_0$, $T_1$, $T_2$, \emph{Urysohn}, $T_3$,
\emph{functionally Hausdorff}, \emph{regular}, \emph{completely regular} and $T_{3.5}$,
that we denote by $\sC_0$, $\sC_1$, $\sC_2$, $\sC_u$, $\sC_3$, $\sC_{fh}$, $\sC_r$, $\sC_{cr}$ and $\sC_{3.5}$ respectively.
\mkp

Take a set $X$ of arbitrary infinite cardinality and an ideal point $p\notin X$ and
set $X^*\defi X\cup \{p\}$ and consider the following two topologies on $X^*$

 $$\tau_{1}=\{U\subseteq X^*: p\notin U\}\cup \{X^*\}$$ $$\tau_{2}=\{U\subseteq X^*: p\in U\}\cup \{\emptyset\}.$$
 \mkp

Set $X_i\defi (X^*,\tau_i)$, $1\leq i\leq 2$ and pick a point  $x$ in $X_1$. Since every  neighbourhood  of $p$ in $X_1$ contains $x$, we have
$\rf_{(X_1,\mathcal{C}_{i})}(x)=\rf_{(X_1,\mathcal{C}_{i})}(p)$ for each
 separation axiom in $\{ T_1$, $T_2$, \emph{Urysohn},
\emph{functionally Hausdorff}, \emph{regular}, \emph{completely regular} and $T_{3.5}\}$. That is to say
 $\rf_{(X_1,\mathcal{C}_{i})}$ is a single-valued map. Now, take $X$, which is a discrete, dense, open subset of $X_1$.
 We have that $\rf_{\mathcal{C}_{i}} X=X\not=\rf_{(X_1,\mathcal{C}_{i})}(X)=\rf_{(X_1,\mathcal{C}_{i})}(p)$.

 As for the $\mathcal{C}_{r}$ (resp. $\mathcal{C}_{cr}$) reflection, remark that no nonempty closed subset of $X_1$ is contained in a proper open subset.
Therefore $\rf_{\mathcal{C}_{r}}X_1=\rf_{\mathcal{C}_{cr}}X_1$ is the indiscrete space
 and again $\rf_{\mathcal{C}_{r}}X=\rf_{\mathcal{C}_{cr}}X=X$, which yields $\rf_{(X_1,\mathcal{C}_{r})}(X)\neq \rf_{\mathcal{C}_{r}}$
and $\rf_{(X_1,\mathcal{C}_{cr})}(X)\neq \rf_{\mathcal{C}_{cr}}$. This completes the proof for open subsets.
For closed subsets, it suffices to take the space $X_2$.

\bdfn  Let $X$ be a topological space and let $\sA$ be a subset of $\topo$. A subset $A$ of $X$ is said
\emph{$\sA$-oset} if there is a space $Y\in \sA$ and a continuous map $f\colon X\longrightarrow Y$ such that
$A=f^{-1}(U)$ for some open subset $U$ of $Y$. It is clear that the family of all $\sA$-osets forms a subbase
in $X$ for the initial topology $\tau_\sA$, with respect to all continuous maps
$X\to Y\in\sA$. The subsets $G\in \tau_\sA$ are called \emph{$\sA$-open}.
A subset $F$ of $X$ is said \emph{$\sA$-closed},  if $X\setminus F$ is $\sA$-open.

In case $\sA=\mathcal{C}_{i}$ for some $i\in \{T_1,T_2,u,r,fh,cr,T_{3.5}\}$,
 we will use the symbolism \emph{$\mathbf T_{i}$-open} for short.
\edfn
\mkp

\blem\label{le_subspace}
Let $X$ be a topological space and let $\sC=\sC(\sA)$ be an epireflective subcategory of $\topo$ that is generated by
a family of spaces $\mathcal A\subseteq \topo$. Given a subset $A$ of $X$, the following assertions are equivalent:
\begin{enumerate}[(a)]
\item $A$ is $\sC$-open (res. $\sC$-closed).
\item $A=\rf_{(X,\sC)}^{-1}(U)$, for some open (resp. closed) subset $U$ of $\rf_\sC X$.
\item $A$ is $\sA$-open (res. $\sA$-closed)
\end{enumerate}
\elem
\bpf
$(a)\Rightarrow (b)$. Let $A$ be a $\sC$-oset subset of $X$. Then there are $Y\in \sC$, $V$ open in $Y$ and a continuous map $f\colon X\longrightarrow Y$
such that $f^{-1}(V)=A$. By the functorial definition of epireflections, there exists a continuous map $g\colon \rf_\sC X\longrightarrow Y$
such that $g\circ \rf_{(X,\sC)}=f$. Therefore $A=f^{-1}(V)=\rf_{(X,\sC)}^{-1}(g^{-1}(V))$ and it
suffices to take $U=g^{-1}(V)$. Now if $A\subseteq X$ is $\sC$-open, then $A =\bigcup\limits_{i\in I}(\bigcap\limits_{j\in J_i}A_{ij})$,
where each $J_i$ is finite and $A_{ij}=\rf_{(X,\sC)}^{-1}(U_{ij})$ for some open subset $U_{ij}$ of $\rf_\sC X$.
Then
$$A =\bigcup\limits_{i\in I}(\bigcap\limits_{j\in J_i}\rf_{(X,\sC)}^{-1}(U_{ij}))=\rf_{(X,\sC)}^{-1}(\bigcup\limits_{i\in I}(\bigcap\limits_{j\in J_i}U_{ij}))$$
and we may choose $$U\defi \bigcup\limits_{i\in I}(\bigcap\limits_{j\in J_i}U_{ij}).$$

$(b)\Rightarrow (c)$. We have seen in Section 2 that the space $\rf_\sC X$ can be realized as the diagonal of a product $\Pi_E X=\prod_{f\in E} Y_f$
where $Y_f\in \sA$ for all $f\in E$ and $f$ stands for a surjective continuous map $f\colon X\to Y_f$. Thus, the family
$\{\pi_f^{-1}(W)\cap \rf_\sC X :\ W\ \hbox{open in}\ Y_f,\ f\in E\}$ form an open subbase in $\rf_\sC X$ and,
as a consequence, the topologies $\tau_\sC$ and $\tau_\sA$ coincide. Therefore, if
$A=\rf_{(X,\sC)}^{-1}(U)$ for some open subset $U$ of $\rf_\sC X$, it follows that $A\in\tau_\sA$.

$(c)\Rightarrow (a)$ is obvious.
\epf
\mkp

The following result characterizes when an epireflection functor preserves subspaces.

\bprp\label{Pr_subspace}
Let $\sC$ be an epireflective subcategory of $\topo$ and let $X$ be a topological space. If $A$ is a subspace of $X$
we have that the epireflection $\sC$ preserves the subspace $A\hookrightarrow X$
if and only if the following two properties are satisfied:
\begin{enumerate}
\item For all $a_1, a_2$ in $A$ such that $\rf_{(A,\sC)}(a_1)\not=\rf_{(A,\sC)}(a_2)$,
we have $\rf_{(X,\sC)}(a_1)\not=\rf_{(X,\sC)}(a_2)$;
\item For every $\sC$-closed (resp. $\sC$-open) subset $F$ of $A$, there is a $\sC$-closed (resp. $\sC$-open) subset $E$ of $X$
such that $E\cap A=F$.
\end{enumerate}
\eprp
\bpf
Suppose that $\sC$ preserves the subspace $A\hookrightarrow X$.
For simplicity, we assume that the $\sC$-image of the subspace inclusion $A\hookrightarrow X$
is the subspace inclusion $\rf_\sC A\hookrightarrow \rf_\sC X$ (that is $\rf_\mathcal{C}A=\rf_{(X,\sC)}(A)$).
Then (1) is obviously satisfied. As for (2) let $U$ be a $\sC$-open.
subset of $A$. By Lemma \ref{le_subspace}, there is an open subset $V$ in $\rf_\sC A$ such that $U=\rf_{(A,\sC)}^{-1}(V)$.
Furthermore, by hypothesis, we may assume that there is an open set $W$ in $\rf_\sC X$ such that
$V=\rf_{\sC}A\cap W$. Thus
$$U=\rf_{(A,\sC)}^{-1}(V)=\rf_{(A,\sC)}^{-1}(\rf_\sC A\cap W)\subseteq \rf_{(X,\sC)}^{-1}(\rf_\sC A)\cap \rf_{(A,\sC)}^{-1}(W)
=\rf_{(X,\sC)}^{-1}(\rf_\sC A)\cap \rf_{(X,\sC)}^{-1}(W).$$ Since also $U\subseteq A$, therefore we have
$$U\subseteq A\cap\rf_{(X,\sC)}^{-1}(\rf_\sC A)\cap\rf_{(X,\sC)}^{-1}(W),$$ and, conversely,
$$U=\rf_{(A,\sC)}^{-1}(V)=\rf_{(A,\sC)}^{-1}(\rf_\sC A\cap W)\supseteq A\cap \rf_{(X,\sC)}^{-1}(W),$$ since $a\in A$ and
$\rf_{(X,\sC)}(a)=\rf_{(A,\sC)}(a)\in W$ imply $\rf_{(A,\sC)}(a)\in \rf_\sC A\cap W$. This proves that
$U=\rf_{(X,\sC)}^{-1}(W)\cap A$.

Conversely, suppose that (1) and (2) are satisfied. Let $j : A\hookrightarrow X$ be the inclusion of the subspace $A$ to
$X$ and let $g\colon \rf_\mathcal{C} A\longrightarrow \rf_{\sC}X$ be a
continuous mapping that makes the following diagram

$$
\xymatrix{ A\ar[dd]_{\rf_{(A,\sC)}} \ar [r]^{j}& X\ar [dd]^{\rf_{(X,\sC)}}\\ \\
\rf_\sC A\ar[r]^{g} & \rf_\sC X}$$
\skp

\noindent commutative

It is clear, by (1), that $g$ is injective. We will show that $g$ is also open.

Let $V$ be an open subset in $\rf_\sC A$. Then  $\rf_{(A,\sC)}^{-1}(V)$ is $\sC$-open in $A$. By (2),
there is a $\sC$-open subset $U$ of $X$ such that $\rf_{(A,\sC)}^{-1}(V)=j^{-1}(U\cap j(A))$.
By the commutativity of the diagram, surjectivity of $\rf_{(A,\sC)}$ and $\rf_{(X,\sC)}$,
and taking into account that $U$ is the inverse image of an open set in $\rf_\sC X$,
it follows that $g(V)=  \rf_{(X,\sC)}(U)\cap \rf_{(X,\sC)}(A)$,
which is open in $\rf_\sC A$. This completes the proof.
\epf

\bdfn
Let $X$ be a topological space and let $\sA$ denote a class of topological spaces, we say that a subset $A$ of $X$ is \emph{$\sA$-embedded}
when for every continuous map $f\colon A\to Z$, with $Z\in\sA$, there exists $Y\in \sA$, with $Z\subseteq Y$, a continuous map $\overline{f}\colon X\to Y$
such that $\overline{f}_{|A}=f$. In other words, every continuous map on $A$ taking values in a space in $\sA$ can be extended
to a continuous map on $X$ taking values in a possibly different larger space in $\sA$.
\edfn

\bcor\label{Co_subspaces}
Let $X$ be a topological space and let $\sC=\sC(\sA)$ be the epireflective subcategory of $\topo$ that is generated by
a $\mathcal A$. If $A$ is an $\sA$-embedded subspace of $X$, then the epireflection functor associated to $\sC$
preserves the subspace $A\hookrightarrow X$.
\ecor
\bpf
We apply Proposition \ref{Pr_subspace} for the proof. In order to verify assertion (1),
suppose that $a_1, a_2$ are two points in $A$ such that $\rf_{(A,\sC)}(a_1)\not=\rf_{(A,\sC)}(a_2)$.
We have seen in Section 2 that the space $\rf_\sC A$ can be realized as the diagonal of a product $\Pi_E X=\prod_{f\in E} Y_f$
where $Y_f\in \sA$ for all $f\in E$ and $f$ stands for a surjective continuous map $f\colon A\to Y_f$. Thus,
there is a map $f\colon A\to Z\in\sA$ such that $f(a_1)\not=f(a_2)$. Then there exists $Y\in \sA$, with $Z\subseteq Y$ and
a continuous map $\overline{f}\colon X\to Y$ that makes following  diagram

$$
\xymatrix{ A\ar[dd]_{j_{(A,X)}} \ar [r]^{f}& Z\ar [dd]^{j_{(Z,Y)}}\\ \\
X\ar [r]^{\overline f} & Y}$$
\skp

\noindent commutative. Thus, we have that $\overline f(a_1)\not=\overline f(a_2)$, which implies that
$\rf_{(X,\sC)}(a_1)\not=\rf_{(X,\sC)}(a_2)$.

As for assertion (2), it suffices to observe that
the collection of $\sA$-osets forms an open subbase for the topology $\tau_\sA =\tau_\sC$.
\epf
\mkp

In \cite{T4}, M. Tkachenko proves that the $\sC_0$-reflection respects arbitrary subgroups. The next corollary
improves this result.

\bcor\label{Ex_T0}
The $\sC_0$-reflection preserves subspaces.
\ecor
\mkp

Even though the $\sC_1$-reflection does not preserve subspaces, our methods provide a neat characterization of this property.
First we need the following lemma.

\blem\label{Le_subspaces}
Let $A$ be a subset of a topological space $X$. Then $A$ is $T_1$-closed if and only if there is a continuous mapping $f\colon X\to Y_{cof}$,
where $Y_{cof}$ is a set equipped with the cofinite topology, such that $A=f^{-1}(p)$ for a singleton $p\in Y$.
 \elem
 \bpf
Suppose that $A$ is $T_1$-closed. Then there is a $T_1$-space $Z$ and a continuous mapping  $g\colon 	X\longrightarrow Z$ such that
$A=g^{-1}(B)$ for a closed subset $B$ in $Z$. If we identify $B$ with a singleton, say $p_B$, and  define $Y\defi (Z\setminus B)\cup \{p_B\}$.
Then the map $f\colon X\to Y_{cof}$, defined by $f(x)=g(x)$ if $x\notin g^{-1}(B)$ and $f(x)=a_B$ if $x\in g^{-1}(B)$, is continuous and
$A=f^{-1}(p_B)$. The converse implication is obvious.
\epf
\mkp

Since the category $\mathbf{Top}_1$ of $T_1$-spaces is generated by the spaces equipped with the cofinite topology, we see that a subset $A$
of a topological space $X$ is \emph{$T_1$-embedded} if and only if for every continuous map $f\colon A\to Z_{cof}$, where $Z$ is a set
equipped with the cofinite topology, there exists a set $Y$, with $Z\subseteq Y$, a continuous map $\overline{f}\colon X\to Y_{cof}$
such that $\overline{f}_{|A}=f$.

\blem\label{Le_subspaces2}
Every $T_1$-closed subspace $A$ of a topological space $X$ is $T_1$-embedded.
\elem
\bpf
Let $f\colon A\to Z^1_{cof}$ be a continuous map defined on $A$.
By Lemma \ref{Le_subspaces}, there is a continuous map $g\colon X\to Z^2_{cof}$ and a point $p\in Z^2$ such that $A=g^{-1}(p)$.
Set $Y\defi Z^1\sqcup Z^2$, the disjoint union of $Z^1$ and $Z^2$, and define $\widetilde{f}\colon X\to Y_{cof}$ by
$\widetilde{f}(x)=f(x)$ if $x\in A$ and $\widetilde{f}(x)=g(x)$ if $x\notin A$. It is clear that inverse image by $\widetilde{f}$ of
every singleton in $Y$ is closed in $X$, which yields the continuity of the map.
\epf

\bthm\label{T1}
If $A$ is a $T_1$-closed subspace of a topological space $X$, then $\rf_\mathcal{C_1}A=\rf_{(X,\mathcal{C_1})}(A).$
\ethm
\bpf
Apply Corollary \ref{Co_subspaces} and Lemma \ref{Le_subspaces2}.
\epf
\mkp

In \cite[Lemma 3.7]{T4} M. Tkachenko proved that the $T_1$-reflection preserves closed subgroups in the category of semitopological groups.
The next Corollary is a variant of Tkachenko's result. Again, our formulation is somewhat more general. First, we recall that a
\emph{left topological group} (resp. \emph{right topological group}) is a group $G$ equipped with a topology such that the left translations
$x\mapsto ax$ are continuous (resp. the right translations $x\mapsto xa$ are continuous).

\bcor
Let $G$ be a left (resp. right) topological group and let $H$ be a closed subgroup of $G$. Then
$\rf_{\mathcal{C_1}}H=\rf_{(G,\mathcal{C_1})}(H).$
\ecor
\bpf
There is no loss of generality in assuming that $G$ is a left topological group.
Let $X\defi (G/H)$  be the quotient space $G/H$. It is easily seen that $X$ is a $T_1$-space.
Furthermore, if $\pi\colon G\to X$ denotes the canonical quotient map,
 we have that $H=\pi^{-1}(\pi(H))$, which implies that $H$ is $T_1$-closed in $G$. Thus it suffices apply Theorem \ref{T1}.
\epf
\mkp

\section{Coincidence of epireflections}

The following general question is dealt with in this section: Let $\sC$ and $\sE$ be two epireflective subcategories of $\topo$
such that $\sC\supseteqq\sE$.
Characterize the spaces $X$ such that $\rf_\sC X=\sE(X)$. This topic has been studied in \cite{Echi:MPRIA,Echi:TP} where it is left as a
specific open question to characterize the spaces $X$ for which $\rf_{\sC_1} X=\rf_{\sC_{3.5}} X$, where $\rf_{\sC_1}$ and $\rf_{\sC_{3.5}}$
are the epireflection functors associated to $\sC_1$ and $\sC_{3.5}$, the subcategories of $T_1$-spaces and $T_{3.5}$ spaces, respectively.

Our approach is based on the notion of $\sC$-open subset that has been introduced previously. First, we recall that
in the category $\topo_0$ the epimorphisms are not the surjections,
like they are for $\topo$. Hence also epireflective has another meaning (cf. \cite{Baron}).

\bthm\label{th_coincidencia}
Let $\sC$ and $\sE$ be two epireflective subcategories of $\topo$
such that $\topo_0\supseteq\sC\supseteq\sE$ and let $X$ be a topological space. Then $\rf_\sC X=\rf_\sE X$ if and only if
every $\sC$-open subset of $X$ is $\sE$-open.
\ethm
\bpf
The ``only if'' part is obvious. As for the ``if'' part,
consider the following commutative diagram

\[
\xymatrix{ & \ar[dl]_{\rf_{(X,\sC)}} X \ar[dr]^{\rf_{_{(X,\mathcal{E})}}} &\\ \rf_\sC X \ar@{>}[rr]_{\rf_{(\rf_\sC \mbox{ X},\sE)}}  & & \rf_\sE X }
\]
\mkp

\noindent where $\rf_{(\rf_\sC { X},\sE)}$ is the unique continuous map canonically defined since $\sC\supseteqq\sE$.  It will suffice to verify that $\rf_{(\rf_\sC { X},\sE)}$ is 1-to-1 and
open. Suppose first that $x,y$ belong to $X$ and $\rf_{(X,\sC)}(x)\not=\rf_{(X,\sC)}(y)$. By our initial assumption  $\sC$ is included in $\topo_0$.
Thus there is an open subset $W$ in $\rf_\sC X$ that contains exactly one of these points. Assume wlog that $\rf_{(X,\sC)}(x)\in W\not\ni \rf_{(X,\sC)}(y)$, which
yields $x\in \rf_{(X,\sC)}^{-1}(W)\not\ni y$. By hypothesis and Lemma \ref{le_subspace},
there must be an open subset $V$ of $\rf_\sE X$ such that
$\rf_{(X,\sC)}^{-1}(W)=\rf_{(X,\sE)}^{-1}(V)$.
Thus $\rf_{(X,\sE)}(x)\not=\rf_{(X,\sE)}(y)$, which proves the injectivity of $\rf_{(\rf_\sC { X},\sE)}$.

Now, let $W$ be an arbitrary open subset of $\rf_\sC X$. Again, there must be an open subset $V$ of $\rf_\sE X$ such that
$\rf_{(X,\sC)}^{-1}(W)=\rf_{(X,\sE)}^{-1}(V)$. Furthermore, the commutativity of the diagram above implies that
$\rf_{(X,\sC)}^{-1}(W)=\rf_{(X,\sC)}^{-1}(\rf_{(\rf_\sC { X},\sE)}^{-1}(V))$. This implies that
$W=\rf_{(\rf_\sC { X},\sE)}^{-1}((V))$ and, as a consequence that $\rf_{(\rf_\sC { X},\sE)}(W)=V$.
This completes the proof.
\epf
\mkp

The following result answers Question 1.6 in \cite{Echi:MPRIA}, repeated in \cite[Question 1.9]{Echi:TP}.

\bcor
Let $X$ be a topological space. Then $\rf_{\sC_1} X=\rf_{\sC_{3.5}} X$ if and only if every $T_1$-closed subset $F$ of $X$
is completely separated from any point $x\notin F$.
\ecor
\bpf
Necessity: Suppose that $\rf_{\sC_1} X=\rf_{\sC_{3.5}} X$. By Theorem \ref{th_coincidencia}, every $T_1$-closed subset
$F$ of $X$ is $\sC_{3.5}$-closed. Therefore, by Lemma \ref{le_subspace}, there is a closed subset $E\subseteq \rf_{\sC_{3.5}} X$ such that $F=\rf_{(X,\sC)}^{-1}(E)$.
Thus, $\rf_{(X,\sC)}(x)\notin E$ for all $x\notin F$. Since  $\rf_{\sC_{3.5}} X$ a $T_{3.5}$ space, this implies that $F$
is completely separated  from any point $x\notin F$.

Sufficiency: Let $U$ be a $T_1$-open subset of $X$, we must verify that $U$ is $\sC_{3.5}$-open in order to
apply Theorem \ref{th_coincidencia}. By hypothesis $X\setminus U$ is completely separated from any point $x\in U$.
Hence $X\setminus U=\rf_{(X,\sC_{3.5})}^{-1}(\rf_{(X,\sC_{3.5})}(X\setminus U))$,
which implies that $X\setminus U$ is $\sC_{3.5}$-closed
and therefore $U$ must be a $\sC_{3.5}$-open subset of $X$, which completes the proof.
\epf
\mkp


\end{document}